\documentclass[12pt]{amsart}
\usepackage{amsfonts}
\usepackage{enumerate}
\usepackage{helvet,courier}
\usepackage{mathptmx,amsmath}
\usepackage{type1cm}
\DeclareMathAlphabet{\mathcal}{OMS}{cmsy}{m}{n}
\DeclareMathAlphabet{\mathbbold}{U}{bbold}{m}{n}  
\usepackage{amsthm}
\usepackage{graphicx}       
\usepackage{multicol}
\usepackage{mathrsfs}
\usepackage{hyperref}
\usepackage{amssymb}
\usepackage{verbatim}
\usepackage{esint}

\usepackage[usenames,dvipsnames]{xcolor}

\theoremstyle{plain}
\newtheorem{thm}{Theorem}[section]
\newtheorem{lm}[thm]{Lemma}

\newtheorem{prop}[thm]{Proposition}

\theoremstyle{remark}
\newtheorem{rmk}{Remark}

\theoremstyle{definition}

\newcommand{\bnu}{\begin{enumerate}}
\newcommand{\enu}{\end{enumerate}}
\newcommand{\bpf}{\begin{proof}}
\newcommand{\epf}{\end{proof}}

\newcommand{\q}{\quad}

\newcommand{\sset}{\subset}

\newcommand{\al}{\alpha}
\newcommand{\be}{\beta}
\newcommand{\ga}{\gamma}
\newcommand{\Ga}{\Gamma}
\newcommand{\om}{\omega}

\newcommand{\ep}{\epsilon}

\newcommand{\si}{\sigma}
\newcommand{\tht}{\theta}

\newcommand{\vp}{\varphi}
\newcommand{\de}{\delta}

\newcommand{\bbz}{\mathbb{Z}}

\newcommand{\bbr}{\mathbb{R}}

\newcommand{\rn}{\mathbb{R}^n}

\newcommand{\bbn}{\mathbb{N}}

\newcommand{\les}{\lesssim}

\newcommand{\f}{\frac}

\newcommand{\nf}{\infty}

\newcommand{\tf}{\tfrac}
\newcommand{\wh}{\widehat}

\newcommand{\supp}{\text{supp }}

\begin{document}

\author{Danqing He}
\address{Danqing He, School of Mathematical Sciences,
Fudan University, People's Republic of China}
\email{hedanqing@fudan.edu.cn}

\author{Kangwei Li}

\address{Kangwei Li, 
Center for Applied Mathematics, Tianjin University, Weijin Road 92, 300072 Tianjin, China}
\email{kli@tju.edu.cn}

\author{Jiqiang Zheng}
\address{Jiqiang Zheng,
  Institute of Applied Physics and Computational Mathematics  and National Key Laboratory of Computational Physics\\ Beijing 100088\\ China}
\email{zheng\_jiqiang@iapcm.ac.cn}

\thanks{
All the authors are supported by  National Key R$\&$D Program of China (No. 2021YFA1002500). In addition, D. He is supported by NNSF of China (No. 12322105), the New Cornerstone Science Foundation, and Natural Science Foundation of Shanghai (No. 22ZR1404900 and 23QA1400300). K. Li was supported by NNSF of China (No. 12222114 and 12001400).  J. Zheng was supported by  NNSF of China (No. 12271051).
}

\title{On pointwise convergence of multilinear  Bochner-Riesz means}
\date{}
\maketitle

\begin{abstract}
We improve the range of indices when the multilinear Bochner-Riesz means converges pointwisely. We obtain this result by establishing the $L^p$ estimates and weighted estimates of $k$-linear maximal Bochner-Riesz operators inductively, which is new when $p<2/k$ in higher dimensions. To prove these estimates, we make use of a variant of Stein's square function and its multilinear generalization.
\end{abstract}



\section{Introduction}

Given $\al\ge 0$, the Bochner-Riesz means is defined by 
$$
B^{\al}(f )(x)=\int_{\bbr^{n}}m^\al(\xi) \wh f(\xi)e^{2\pi ix\cdot\xi} d\xi,
$$
where $m^\al(\xi)=(1-|\xi|^2)^\al_+$.
The  maximal Bochner-Riesz operator is defined as
$$
B^{\al}_*(f)(x)=\sup_{R>0}\Big|B^\al_R(f)(x)\Big|,
$$
where
$$
B^{\al}_R(f )(x)=\int_{\bbr^{n}}m^\al(R^{-1}\xi) \wh f(\xi)e^{2\pi ix\cdot\xi} d\xi.
$$
The famous Bochner-Riesz conjecture says that $B^\al(f)$ is bounded on $L^p(\rn) $ if and only if  $\al>\max\{n|\tf1p-\tf12|-\tf12,0\}$. This conjecture is closely related to many other conjectures including the local smoothing conjecture and the restriction conjecture; see \cite{Tao1999a} for details. The boundedness of $B^\al$ and $B^\al_*$ are studied extensively in recent years. One may consult \cite{Bourgain2011, Guo2022, Li2024}, and references therein.

The multilinear theory in harmonic analysis goes back to \cite{Coifman1975}. After the seminal work of Lacey and Thiele \cite{Lacey1997, Lacey1999}, who obtained the boundedness of bilinear Hilbert transform, the multilinear generalization of various operators were studied in past decades, among which are
multilinear singular integrals \cite{Kenig1999, Grafakos2002, Grafakos2015}, multilinear H\"ormander multipliers \cite{Tomita2010}, bilinear Hilbert transform along curves \cite{Li2013, Li2013a}, bilinear square functions \cite{Bernicot2011}, and  bilinear Bochner-Riesz operators \cite{Bernicot2015a, Liu-Wang2020}.

Given $\al\ge0$,
the bilinear Bochner-Riesz operator is defined by
$$
B^{2,\al}(f,g )(x)=\int_{\bbr^{2n}}m^\al(\xi,\eta) \wh f(\xi)\wh g(\eta)e^{2\pi ix\cdot(\xi+\eta)} d\xi\ d\eta,
$$
where $m^\al(\xi,\eta)=(1-|(\xi,\eta)|^2)^\al_+$.
When $\al=0$, this operator was studied by \cite{Grafakos2006a} and \cite{Diestel2007}.
The systematical study of bilinear Bochner-Riesz means when $\al>0$ was initiated by Bernicot, Grafakos, Song and Yan \cite{Bernicot2015a}. 
They established the boundedness of $B^\al$ by estimating the Sobolev norm of the bilinear multiplier $m^\al(\xi,\eta)$ using that it is bi-radial.
For example, they showed that
$$
\|B^{2,\al}(f,g)\|_{L^{n/(n+1)}}\le C\|f\|_{L^{2n/(n+1)}}\|g\|_{L^{2n/(n+1)}}
$$
when $\al>\tf{3n-1}{n+1}$.
Jeong, Lee and Vargas \cite{Jeong2017} improved the boundedness of the bilinear Bochner-Riesz operator by a decomposition lemma, which follows from decomposing $m^\alpha(\xi,\eta)$ along its singularities $\{(\xi,\eta)\in\bbr^{2n}:\ |(\xi,\eta)|=1\}$ and exploring the interaction between $\xi$ and $\eta$. This lemma allows them to obtain the relation between the bilinear Bochner-Riesz means and the classical Stein's square function
$$
G^\al(f)(x)=\Big(\int_0^\nf|(\si^\al(t^{-1}\cdot)\wh f)^\vee (x)|^2\f{dt}t\Big)^{1/2},
$$
where $\si^\al(\xi)=|\xi|^2(1-|\xi|^2)_+^\al$.

\medskip

Similar to the linear case, to study the pointwise convergence of 
 bilinear  Bochner-Riesz means, Grafakos, Honz\'ik and the first author \cite{Grafakos2018a} studied the bilinear maximal Bochner-Riesz operator defined by
$$
B^{2,\al}_*(f, g)(x)=\sup_{R>0}\Big|\int_{\bbr^{2n}}m^\al(R^{-1}\xi,R^{-1}\eta) \wh f(\xi)\wh g(\eta)e^{2\pi ix\cdot(\xi+\eta)} d\xi\ d\eta\Big|.
$$ 
They used  the wavelet decomposition to show that
$B^{2,\al}_*$ is bounded from $L^2(\rn)\times L^2(\rn)$ to $L^1(\rn)$ when $\al>(2n+3)/4$.
Later,  \cite{Jeong2019} improved this result
by developing the method introduced in 
\cite{Jeong2017}. 

Recently, \cite{Kaur2020} introduced a new method to study the boundedness of bilinear maximal Bochner-Riesz means. More precisely, they decomposed $m^{\al}(\xi,\eta)$ along the singularities in $\xi$ first, by which they were able to apply a reproducing formula (Proposition~\ref{tIP1}) to establish a new connection between the bilinear maximal Bochner-Riesz operator and 
Stein's square functions. This method was utilized by  \cite{Shuin2024} to obtain the boundedness of multilinear Bochner-Riesz operators.

\medskip

In this note, we consider the pointwise convergence of multilinear Bochner-Riesz means, which motivated the study of multilinear maximal Bochner-Riesz operators. 

Given $\al\ge 0$, we define the $k$-linear Bochner-Riesz multiplier by
$$
m^\al(\xi)=(1-|\xi|^2)^\al_+
$$
for $\xi=(\xi_1,\xi_2,\dots,\xi_k)\in (\bbr^n)^k$. The associated multilinear Bochner-Riesz means is defined as
$$
B^{k,\al}_R(f_1,\dots,f_k)(x)=\int_{\bbr^{nk}}m^\al(R^{-1}\xi)\prod_{\iota=1}^k \wh f_\iota(\xi_\iota)e^{2\pi ix\cdot(\xi_1+\cdots+\xi_k)} d\xi,
$$
where $x\in\rn$.
The multilinear maximal Bochner-Riesz operator is defined as
$$
B^{k,\al}_*(f_1,\dots, f_k)(x)=\sup_{R>0}\Big|B^{k,\al}_R(f_1,\dots,f_k)(x)\Big|.
$$

Our main result is as follows.

\begin{thm}\label{Main}Let $n\ge 2$. 
		Suppose that $p_\iota\in (1,\nf)$ for $\iota=1,\dots,k$, and 
\begin{equation}\label{ePCBB13}
		\al>\al(p_1,\dots, p_k):=\sum_{\iota=1}^k
		\max\big\{n(\tf1{p_\iota}-\tf12),n(\tf12-\tf1{p_\iota})-\tf12,0\big\}+\tf{k-2}2.
\end{equation}
		Then for $f_\iota\in L^{p_\iota}(\rn)$, $\iota=1,\dots,k$, 
	\begin{equation*} 
		\lim_{R\to\nf}B^{k,\al}_R(f_1,\dots, f_k)(x)= \prod_{\iota=1}^kf_\iota(x)\q a.e. \quad  x\in\rn.
	\end{equation*}
\end{thm}

Our result is new when $p_\iota<2$ for some $\iota$. For example, for $(p_1,p_2)=(1/(1-\ep),2)$ in the bilinear setting, we need only $\al>\tf{n(1-2\ep)}2$, while the known trivial sufficient condition is $\al>\tf{2n-1}2$.

When $p_\iota\ge 2$ for all $\iota$, the range of $\al$ in \eqref{ePCBB13} is larger than the known range of $\al$ when the multilinear maximal Bochner-Riesz is bounded.

To prove Theorem~\ref{Main} when some $p_\iota<2$, we make use of 
$$
\tilde G^{\al}(f)(x)=\sup_{R>0}\Big(\f 1R\int_0^R |B^{\al}_t(f)(x)|^2 dt\Big)^{1/2}
$$ 
other than 
$
G^\al(f).
$
One advantage of $\tilde G^{\al}$ lies in the compactness of $t$ in its definition, which allows us to obtain the boundedness of 
$B^{k,\al}_*$ when $p_\iota<2$. 
Another advantage of $\tilde G^{\al}$ and its $k$-linear generalization comes from the neat iterated relation \eqref{eMBR5} below, which simplifies the inductive argument in the study of multilinear maximal Bochner-Riesz operators.
For the case $p_\iota\ge 2$, we consider the weighted boundedness of $B^{k,\al}_*$ instead of the $L^p$ boundedness, which follows the strategy introduced in \cite{Carbery1988}. More precisely, we have the following result. We denote $I_k=\{1,\dots,k\}$, and  assume that $I\sset I_k$.

\begin{thm}\label{tMBR6}
	Let $\ga_\iota\in(0,n)$  for $\iota\in I$, where $n\ge 2$, 
	 $p_\iota\in (1,2]$ for $\iota\in I_k\setminus I$.
	Suppose 
	 $\tf1p=\sum_{\iota\in I}\tf12+\sum_{\iota\in I_k\setminus I}\tf1{p_\iota}$,
	 $\ga=\tf p{2}\sum_{\iota\in I}\ga_\iota$, and
	$$
	\al>\sum_{\iota\in I}\max\{\tf{\ga_\iota-1}2,0\}+\sum_{\iota\in I_k\setminus I}n(\tf1{p_\iota}-\tf12)+\tf{k-2}2.
	$$
	Then 
	$$
	\| B_*^{k, \alpha}(f_1,\dots, f_k)\|_{L^p(|x|^{-\ga})}
	\les \prod_{\iota\in I}\|f_\iota\|_{L^2(|x|^{-\ga_\iota})}
	\prod_{\iota\in I_k\setminus I} \|f_\iota\|_{L^{p_\iota}(\rn)}.
	$$
\end{thm}

As 
$$
L^q(\rn)\sset L^2(\rn)+L^2(|x|^{-\om})
$$
for $\om>n(1-\tf 2q)$ when $q\ge 2$, one can obtain Theorem~\ref{Main} from Theorem~\ref{tMBR6} directly.

\medskip

Besides the new boundedness of $B^{k,\al}_*$, we have also a counterexample which is useful when $p_\iota$ are small.

\begin{thm}\label{thm:countexa}
Let $n\ge 2$. 
The multilinear maximal estimate
\begin{equation*} 
\big\|B^{k,\al}_*(f_1,\dots, f_k)(x)\big\|_{L^p(\bbr^n)}\leq C\prod_{\iota=1}^k\|f_\iota\|_{L^{p_\iota}(\bbr^n)}
\end{equation*}
fails when
\begin{equation}\label{equ:condialp}
\alpha<\tf{2n-1}{2p}- \tfrac{(n-1)k+1}2.
\end{equation}

\end{thm}

\begin{rmk}

(i)
When $p_\iota$ is close to $1$ for all $\iota$, the gap between the critical index in \eqref{equ:condialp} and the positive result in \eqref{ePCBB13} is 
$$
 \big(\sum_{\iota=1}^k n(\tf1{p_\iota}-\tf12) +\tf{k-2}2\big)-(\tf{2n-1}{2p}- \tfrac{(n-1)k+1}2)=\tf1{2p}-\tf12.
$$

(ii)
When $p>\tf{2(2n-1)}{(n-1)k+1}$, the critical index in \eqref{equ:condialp} is negative, so the counterexample is useless in this case.

\end{rmk}

\medskip

This paper is organized as follows. In Section 2, we discuss some results on Stein's square function $\tilde G^{\al}$ and the strategy in \cite{Kaur2020}. In Section 3, we prove the bilinear version of Theorem~\ref{Main}. In Section 4, we discuss briefly how to obtain Theorem~\ref{Main} inductively. In Section 5, we discuss the counterexample.

\section{Preliminaries}

\subsection{Stein's square function}

Let $\si^\al(\xi)=|\xi|^2(1-|\xi|^2)_+^\al$. The following  two versions of Stein's square functions are closely related to the maximal Bochner-Riesz operator:
$$
\tilde G^{\al}(f)(x)=\sup_{R>0}\Big(\f 1R\int_0^R |B^{\al}_t(f)(x)|^2 dt\Big)^{1/2},
$$ 
$$
G^\al(f)(x)=\Big(\int_0^\nf|(\si^\al(t^{-1}\cdot)\wh f)^\vee (x)|^2\f{dt}t\Big)^{1/2}.
$$

\begin{prop}[{\cite[Section VII.5]{Stein1971a}}]
Let $\al_0>0$, and $1<p<\nf$.
	Suppose
	$$
	\|G^\al\|_{L^p(\rn)\to L^p(\rn)}<\nf\q\forall \al>\al_0.
	$$
	Then
	
	(i) for $\al>\al_0$ we have
	$$
	\|\tilde G^\al\|_{L^p(\rn)\to L^p(\rn)}<\nf;
	$$
	
	(ii) for $\be >\al_0+\tf12$, we have
	$$
	\|B^\be_*\|_{L^p(\rn)\to L^p(\rn)}<\nf.
	$$
	
\end{prop}

The statement (i) follows from the  relation
$$
\tilde G^{\al}(f)\le \sum_{s=0}^N G^{\al+s}(f)+B^{\al+N+1}_*(f),
$$
where $N$ is sufficiently large such that
$B^{\al+N+1}_*(f)\les M_{HL}(f)$ with $M_{HL}$ being the Hardy-Littlewood maximal function. 
(ii) is a consequence of the pointwise control
$$
B^{\al+\de}_*(f)(x)\les \tilde G^\al (f)(x)
$$
for $\de>\tf12$. We refer to \cite[Section VII.5]{Stein1971a} for more details.

It is well known that, when $1< p\le2$,
$$
\|G^\al(f)\|_{L^p(\rn)}\le C\|f\|_{L^p(\rn)}
$$
if and only if $\al>n(\tf1p-\tf12)-\tf12$; see, for example, \cite{Sunouchi1967} and \cite{Igari1971}.

When $p> 2$, it is conjectured that 
\begin{equation*} 
\|G^\al(f)\|_{L^p(\rn)}\le C\|f\|_{L^p(\rn)}\q\forall \al>\al(p)-\tf12,
\end{equation*}
 where 
$\al(p)=\max\{n|\tf1p-\tf12|-\tf12,0\}$. In particular, $\al(p)=\tf{n-1}2-\tf np$ when $p\ge \tf{2n}{n-1}$.

To study the boundedness of $G^\al$, we decompose it dyadically.
We take
$\psi\in\mathcal C_c^\nf(\bbr)$ such that
$\supp\psi\sset(\tf12,2)$ and 
\begin{equation}\label{e000}
	\sum_{j\in\bbz}\psi(2^{-j}\xi)=1\q\forall \xi\neq0.
\end{equation}
For $j\ge 1$, we define
$$
\si^\al_j(\xi)=\si^\al(\xi)\psi(2^j(1-|\xi|^2)),
$$
and
$$
\si^\al_0=\si^\al-\sum_{j\ge 1}\si^\al_j.
$$
Recalling \eqref{e000},
$\si^\al_0$ is of the form
$
\si^\al\vp
$
with 
\[
\vp(\xi)=\sum_{j\le 0}\psi(2^j(1-|\xi|^2))= \psi( 1-|\xi|^2)\in C_c^\infty(\mathbb R^n),
\]where obviously $\supp\vp\sset B(0,\tf1{\sqrt 2})$. 
Correspondingly 
\begin{equation}\label{eLBR2}
G^\al(f)(x)\le\sum_{j\ge 0}2^{-j\al}\Big(\int_0^\nf|(2^{j\al}\si^\al_j(t^{-1}\cdot)\wh f)^\vee (x)|^2\f{dt}t\Big)^{1/2}:=\sum_{j\ge 0}2^{-j\al}G_j^\al(f)(x)  ,
\end{equation}
for $j\ge 0$.

The study of the case $j=0$ goes to the classical Littlewood-Paley theory. 
For $j\ge 1$ we can write 
\begin{align*}
2^{j\al}\si^\al_j(\xi)&= |\xi|^2 2^{j\al}(1-|\xi|^2)_+^\alpha \psi(2^j(1-|\xi|^2))\\
&= |\xi|^2 \psi_1 (2^j(1-|\xi|^2))=  \psi_1 (2^j(1-|\xi|^2))- 2^{-j} \psi_2 (2^j(1-|\xi|^2)) \\
&=:  \widetilde \psi(2^j(1-|\xi|^2)) ,
\end{align*}
where $\psi_1(t)=t^\alpha \psi(t) $, $\psi_2(t)= t \psi_1(t)$ and $\widetilde \psi= \psi_1-2^{-j}\psi_2$.
Since $\widetilde \psi \in C^\nf_c(\bbr)$ and $\supp \widetilde \psi\sset(\tf12,2) $ (in below we may abuse of notation and still use $\psi$ instead of $\widetilde \psi$),
the study of $G_j^\al$ is reduced to 
$$
G_j^\psi(f)(x):=\Big(\int_0^\nf\Big|B^\psi_{2^{-j},t}(f)(x)\Big|^2\f{dt}t\Big)^{1/2},
$$
where 
$$
B^\psi_{2^{-j},t}(f)(x):=[\psi(2^j(1-|\xi/t|^2))\wh f(\xi)]^\vee(x).
$$
Similarly we define also
$$
\tilde G_j^{\psi}(f)(x)=\sup_{R>0}\Big(\f 1R\int_0^R |B^{\psi}_{2^{-j},t}(f)(x)|^2 dt\Big)^{1/2}.
$$

It follows from the classical Littlewood-Paley theory that
$G^\psi_j$ is bounded on $L^p(\rn)$ for $j\ge 1$.

\begin{prop}\label{t001}
	Let $1<p<\nf$, $j\ge 1$. Then there exists a finite constant $C_j$ such that 
	$$
	\|G^\psi_j\|_{L^p(\rn)\to L^p(\rn)}<C_j.
	$$
\end{prop}

  In particular cases we can get quantitative estimate of Proposition \ref{t001}. For instance, 
it is direct to check that 
\begin{equation}\label{eLBR18}
\|G_j^\psi(f)\|_2\le C2^{-j/2}\|f\|_2,
\end{equation}
which by \eqref{eLBR2} implies further that
\begin{equation*} 
\|G^\al(f)\|_{L^2(\rn)}\le C\|f\|_{L^2(\rn)}\q\forall\al>-\tf12.
\end{equation*}

To study the boundedness of $\tilde G^\psi_j$ on $L^p$ for $1<p\le 2$, we need the following result. This result should be known, but we cannot find an exact reference, so we sketch the proof below.

\begin{lm}

Let $K_j(x)=[\psi(2^j(1-|\xi|^2))]^\vee(x)$. Then
\begin{equation}\label{eLBR11}
	|K_j(x)|\le \left\{\begin{array}{ll}
		C_{\psi}2^{-j}(1+|x|)^{-(n-1)/2},&\q |x|\le 2^j\\
		C_{\psi}2^{-j(n+1)/2}(1+2^{-j}|x|)^{-M},&\q |x|\ge 2^j.
	\end{array}\right.
\end{equation}

\end{lm}

\bpf

 Let $\de=2^{-j}$, then
$$
K_j(x)=\int_0^{+\infty} \psi(\de^{-1}-\de^{-1}r^2)\int_{\mathbb S^{n-1}}e^{2\pi irx\cdot\tht}d\si(\tht) r^{n-1}dr.
$$
It follows from the Fourier transform of the unit sphere that, for $|x|\ge1$, the main term of $K_j$ is
$$
\int_0^{+\infty}  \psi(\de^{-1}-\de^{-1}r^2) |rx|^{-(n-1)/2}e^{2\pi i|rx|} r^{n-1}dr.
$$
Applying integration by parts, last expression is bounded
by 
$|x|^{-\tf{n-1}2}(\de|x|)^{-M}\de$ for any $M\ge 0$, since
$\supp\psi(\de^{-1}-\de^{-1}r^2)\sset [1-3\de,1]$.
In particular, for $1\le |x|\le\de^{-1}$, we have
$$
|K_j(x)|\le C2^{-j}|x|^{-(n-1)/2}
$$
by taking $M=0$; for $|x|\ge\de^{-1}$, we have 
$$
|K_j(x)|\le C\de^{\tf{n+1}2}(1+\de|x|)^{-M}.
$$
When $|x|\le1$, we have the trivial bound 
$$
|K_j(x)|\le C2^{-j}
$$
as $\text{supp } K_j\sset S^{n-1}+O(2^{-j})$.

\epf

The following result is important in studying the $L^{p_1}\times\cdots\times L^{p_k}\to L^p$ boundedness of $B^{k,\al}_*$ when some $p_\iota<2$.

\begin{prop}
Let $1< p\le 2$. For any $\ep>0$, we can find $C_\ep>0$ such that

(i)
\begin{equation}\label{eLBR14}
\|\sup_{t>0}|B_{2^{-j},t}^\psi(f)|\|_{L^p(\rn)}\le C_\ep 2^{j((n-1)(\tf1p-\tf12)+\ep)}\|f\|_{L^p(\rn)};
\end{equation}

(ii) 
\begin{equation}\label{eLBR19}
	\|\tilde G_j^{\psi}(f)\|_{L^p(\rn)}
\le C_\ep 2^{j(n(\tf1p-\tf12)-\tf12+\ep)}\|f\|_{L^p(\rn)}.
\end{equation}

(iii) For $\de>n(\tf1p-\tf12)-\tf12$, we have
\begin{equation}\label{eLBR20}
\|\tilde G^\de(f)\|_{L^p(\rn)}\le C\|f\|_{L^p(\rn)}.
\end{equation}
\end{prop}

\bpf

(i)
As the right hand side of \eqref{eLBR11} is a radial decreasing function whose integral is $C_\psi 2^{j(n-1)/2}$, we obtain from 
\cite[Theorem 2.1.10]{Grafakos2014b} that
$$
|K_{j}*f(x)|\le C_\psi 2^{j(n-1)/2} M_{HL}(f)(x).
$$
As a result, for any $p>1$, we have
\begin{align*}
	\Big\|\sup_{t>0}|B_{2^{-j},t}^\psi(f)|\Big\|_{L^p(\rn)}
	\le C2^{j(n-1)/2} \|M_{HL}(f)\|_{L^p(\rn)}
	\le  C2^{j(n-1)/2} \|f\|_{L^p(\rn)}.
\end{align*}
For $p=2$,  note that we can write 
\begin{align*}
(B_{2^{-j},t}^\psi(f))^2= \int_0^t 2 B_{2^{-j},u}^\psi(f) \frac{d}{du}(B_{2^{-j},u}^\psi(f)) du.
\end{align*}
Direct computations give that 
\begin{align*}
\frac{d}{du}(B_{2^{-j},u}^\psi(f))=2^{j+1}u^{-1} [\widetilde\psi(2^j(1-|\xi/u|^2))\wh f(\xi)]^\vee(x),
\end{align*}
where 
\[
\widetilde\psi(t)= -2^{-j} t \psi'(t)+ \psi'(t).
\]
Thus by H\"older's inequality, we see that 
\[
|B_{2^{-j},t}^\psi(f)|^2\lesssim 2^j  G_j^{\psi} (f) G_j^{\widetilde\psi}(f).
\]
Then 
it follows from \eqref{eLBR18} that
$$
\|\sup_{t>0}|B_{2^{-j},t}^\psi(f)|\|_{L^2(\rn)}\les \|f\|_{L^2(\rn)}.
$$
Interpolating between above two estimates yields the conclusion for $p\in(1,2)$.

(ii) For any $p\in (1,2)$, we apply \eqref{eLBR14}
to obtain 
$$
\|\tilde G_j^{\psi}(f)\|_{L^p(\rn)}
\le C_\ep 2^{j((n-1)(\tf1p-\tf12)+\ep)}\|f\|_{L^p(\rn)}
$$
since $\tilde G_j^{\psi}(f)\les \sup_{t>0}|B_{2^{-j},t}^\psi(f)(x)|$.
For $p=2$, we may apply \eqref{eLBR18} to obtain 
$$
\|\tilde G_j^{\psi}(f)\|_{L^2(\rn)}
\le C 2^{-\tf j2}\|f\|_{L^2(\rn)},
$$
as 
$$
\tilde G_j^{\psi}(f)(x)\le G_j^{\psi}(f)(x).
$$
Interpolating between above two estimates yields the conclusion for $p\in(1,2)$.

(iii) follows from (ii) and 
$$
\tilde G^\de(f)(x)\le \sum_{j\ge 0}2^{-j\de}\tilde G^\psi_j(f)(x)
$$
for appropriate $\tilde G^\psi_j$.

\epf

The following  weighted estimate is essentially proved in \cite{Carbery1988}.

\begin{thm}\label{tLBR5}
	Suppose $j\in\bbn$, and $\ga\in[0,n)$. Then
	\begin{equation}\label{eLBR16}
	\|G_j^\psi (f)\|_{L^2(|x|^{-\ga})}\le CA_j(\ga)\|f\|_{L^2(|x|^{-\ga})},
	\end{equation}
	where
	$$
	A_j(\ga)=\left\{\begin{array}{ll}
		2^{-j(2-\ga)/2},&\q\ga\in(1,n)\\
		(j2^{-j})^{1/2},&\q\ga=1\\
		2^{-j/2},&\q\ga\in[0,1).
	\end{array}\right.
	$$
	As a result,
\begin{equation}\label{eLBR8}
		\|G^\de(f)\|_{L^2(|x|^{-\ga})}\le C\|f\|_{L^2(|x|^{-\ga})}
\end{equation}
	for $\de>\max\{\tf{\ga-1} 2,0\}-\tf12$.
\end{thm}

\begin{rmk}
	We may rewrite \eqref{eLBR16} as 
	$$
	\|G_j^\psi (f)\|_{L^2(|x|^{-\ga})}\le C_\ep 2^{j(\max\{\tf{\ga-1}2,0\}-\tf12+\ep)}\|f\|_{L^2(|x|^{-\ga})}.
	$$
\end{rmk}


\subsection{Bilinear maximal Bochner-Riesz and bilinear square functions}
The bilinear maximal Bochner-Riesz operators and related bilinear square functions are studied by \cite{Kaur2020} and \cite{Choudhary2023}.

Let us review the definitions and their strategies. We will follow their strategy to obtain pointwise convergence of bilinear Bochner-Riesz means via establishing weighted and unweighted estimates. 

Given $\al>0$, the bilinear maximal Bochner-Riesz operator is defined by 
$$
B^{2,\al}_*(f_1,f_2)(x)=\sup_{R>0}
\Big|\int_{\bbr^{2n}}m^{\al}(R^{-1}\xi)\wh{f_1}(\xi_1)\wh{f_2}(\xi_2)e^{2\pi ix\cdot(\xi_1+\xi_2)}d\xi
\Big|.
$$

We can decompose $m^{\al}$ into
$$
m^\al_0(\xi)+\sum_{j\ge 1}m^\al_j(\xi),
$$
where
$$
m^\al_j(\xi)=(1-|\xi|^2)^{\al}_+\psi(2^j(1-|\xi_2|^2)),
$$
and 
\begin{equation}\label{eLBR5}
m^\al_0=m^\al-\sum_{j\ge 1}m^\al_j.
\end{equation}

We then need to estimate the $L^p$-norm of 
$$
M_j(f_1,f_2):=\sup_{R>0}
\Big|\int_{\bbr^{2n}}m^\al_j(R^{-1}\xi)\wh{f_1}(\xi_1)\wh{f_2}(\xi_2)e^{2\pi ix\cdot(\xi_1+\xi_2)}d\xi\Big|.
$$

The following reproducing formula plays an important role in the study of bilinear maximal Bochner-Riesz operators.

\begin{prop}[{\cite[p 278]{Stein1971a}}]\label{tIP1}
	Let $\be>0$ and $\de>-1$. Then
	$$
	(1-\tf{|\eta|^2}{R^2})_+^{\be+\de}=
	C_{\be,\de}R^{-2\be-2\de}\int_{|\eta|}^R (R^2-t^2)^{\be-1}t^{2\de+1}(1-\tf{|\eta|^2}{t^2})^\de dt,
	$$
	where $C_{\be,\de}=\tf{2\Ga(\de+\be+1)}{\Ga(\de+1)\Ga(\be)}$.
\end{prop}

Let $\al=\de_1+\be_1>0$, $\de_1>-\tf12$,  $\be_1>0$.
Write 
\[
\left( 1- \frac{|\xi_1|^2+|\xi_2|^2}{R^2}\right)_+^\alpha= \left( 1-\frac{|\xi_2|^2}{R^2}\right)_+^\alpha \left( 1-\frac{|\xi_1|^2}{R^2-|\xi_2|^2}\right)_+^\alpha. 
\]
Applying Proposition~\ref{tIP1} with  $(|\eta|,R)$ replaced by $(|\xi_1|, \sqrt{R^2-|\xi_2|^2})$ we get that
\begin{align}
\Big|\int_{\bbr^{2n}}m^\al_j(R^{-1}\xi)&\wh{f_1}(\xi_1)\wh{f_2}(\xi_2)e^{2\pi ix\cdot(\xi_1+\xi_2)}d\xi\Big|\notag\\
\les&2^{-j/4}\tilde G^{\de_1}(f_1)(x)\Big(\int_0^{2^{-(j-1)/2}}|B_{j,t,R}^{\be_1}(f_2)(x)t^{2\de_1+1}|^2 dt\Big)^{1/2},\label{eLBR21}
\end{align}
where 
\begin{equation*} 
B^{\be}_{j,t,R}(f)(x)=\int_{\rn}(1-t^2-\tf{|\eta|^2}{R^2})^{\be-1}_+ \psi(2^j(1-\tf{|\eta|^2}{R^2}))\wh f(\eta)e^{2\pi ix\cdot\eta}d\eta.
\end{equation*}
In particular, 
\begin{equation}\label{eLBR7}
	M_j(f_1,f_2)(x)\le 2^{-j/4}\tilde G^{\de_1}(f_1)(x)B^{\be_1,\de_1}_{j,*}(f_2)(x),
\end{equation}
where 
\begin{equation}\label{eLBR6}
B^{\be_1,\de_1}_{j,*}(f_2)(x)=\sup_{R>0}\Big(\int_0^{2^{-(j-1)/2}}|B_{j,t,R}^{\be_1}(f_2)(x)t^{2\de_1+1}|^2 dt\Big)^{1/2}.
\end{equation}

The boundedness of $B^{\be_1,\de_1}_{j,*}$ is closely related to the boundedness of $G_j^\psi$. To state it accurately, we define 
$$
\al_{p_c}(p):=\al(p_c)\f{2^{-1}-p^{-1}}{2^{-1}-p_c^{-1}}
$$
for $p_c\ge\tf{2n}{n-1}$ and $2\le p\le p_c$, where recall that $\al(p_c)=n(\tf12-\tf1{p_c})-\tf12$.

\begin{prop}[{\cite[Theorem 5.1]{Kaur2020}}]\label{tLBR2}
Let $p_c>\tf{2n}{n-1}$.
Suppose that for any $\ep>0$, there exists some $C_\ep>0$ such that
\begin{equation}\label{eLBR3'}
\Big\|\Big(\int_1^2\Big|[\psi(2^j(t^2-|\xi|^2))\wh f(\xi)]^\vee(x)\Big|^2dt\Big)^{1/2}\Big\|_{L^p(\rn)}\le C_\ep 2^{-j(\tf12-\al(p)-\ep)}\|f\|_{L^p(\rn)}
\end{equation}
holds for all $p>p_c$.
Then for $p\in(p_c,\nf)$ and $\be>\al(p)+\tf12$, we have
$$
\|B^{\be,\de}_{j,*}(f)\|_{L^p(\rn)}\le C_\ep 2^{-j(\al-\al(p)-\tf14-\ep)}\|f\|_{L^p(\rn)},
$$
where $\al=\be+\de$;
for $p\in(2,p_c)$ and $\be>\al_{p_c}(p)+\tf12$, we have
$$
\|B^{\be,\de}_{j,*}(f)\|_{L^p(\rn)}\le C_\ep 2^{-j(\al-\al_{p_c}(p)-\tf14-\ep)}\|f\|_{L^p(\rn)}.
$$

\end{prop}

Assuming \eqref{eLBR3'},
we can now apply \eqref{eLBR7} and Proposition~\ref{tLBR2} to control
$
\|M_j\|_{L^{p_1}\times L^{p_2}\to L^p},
$
which in turn leads to the control of 
$$
\|B^{2,\al}_*\|_{L^{p_1}\times L^{p_2}\to L^p}.
$$

Let 
$$
\al_{p_c}(p_1,p_2)=
\left\{
\begin{array}{ll} 
\al(p_1)+\al(p_2)&p_c\le p_1,p_2<\nf\\
\al(p_1)+\al_{p_c}(p_2)&p_1\in[p_c,\nf),\ p_2\in[2,p_c)\\
\al_{p_c}(p_1)+\al(p_2)&p_1\in[2,p_c),\ p_2\in[p_c,\nf)\\
\al_{p_c}(p_1)+\al_{p_c}(p_2)&p_1,p_2\in[2,p_c).
\end{array}
\right.$$

\begin{thm}[{\cite[Theorem 2.1]{Kaur2020}}]
Let $p_c>\tf{2n}{n-1}$.
Suppose \eqref{eLBR3'} holds for all $p>p_c$.
Then, for any $p_1,p_2\in[2,\nf)$, we have
$$
\|B_*^{2,\al}(f_1,f_2)\|_p\le C\|f_1\|_{p_1}\|f_2\|_{p_2}
$$
whenever $\al>\al_{p_c}(p_1,p_2)$.

\end{thm}


\section{Pointwise convergence of bilinear Bochner-Riesz means}

In this section we study the pointwise convergence of $B^{2,\al}_R$ as $R\to\nf$. 
For $1<p_1,\ p_2<\nf$, and $\tf1p=\tf1{p_1}+\tf1{p_2}$, we define
\begin{align*}
\al(p_1,p_2)=&\sum_{\iota=1}^2
\max\{n(\tf1{p_\iota}-\tf12),n(\tf12-\tf1{p_\iota})-\tf12,0\}\\
=&
\left\{
\begin{array}{ll} 
\sum_{\iota=1}^2\max\{n(\tf12-\tf1{p_\iota})-\tf12,0\},&2\le p_1,p_2<\nf\\
n(\tf1{p_2}-\tf12)+\max\{0,n(\tf12-\tf1{p_1})-\tf12\}, &1< p_2\le2\le  p_1<\nf\\
n(\tf1{p_1}-\tf12)+\max\{0,n(\tf12-\tf1{p_2})-\tf12\}, &1< p_1\le2\le  p_2<\nf\\
n(\tf1p-1), &p_1,p_2\in(1,2].\\
\end{array}
\right.
\end{align*}

Our main result is Theorem~\ref{Main} when $k=2$. Let us state it in the bilinear setting.

\begin{thm}\label{tPCBB1}
	Let $1<p_1,p_2<\nf$, and $\al>\al(p_1,p_2)$. Suppose $f_\iota\in L^{p_\iota}(\rn)$ for $\iota=1,2$, where $n\ge 2$. Then
	\begin{equation}\label{ePCBB10}
	B^{2,\al}_R(f_1,f_2)(x)\to f_1(x)f_2(x)\q a.e.
	\end{equation}
\end{thm}

\begin{rmk}
In the linear setting, it is known that, if 
$\al>\max\{(n-1)(\tf1p-\tf12), n(\tf12-\tf1p)-\tf12,0\}$, then
$B^\al(f)\to f$ a.e. for $f\in L^p(\rn)$. However, in the bilinear setting, one cannot expect that a similar condition
$$ 
\al>\sum_{j=1}^2\max\{(n-1)(\tf1{p_j}-\tf12), n(\tf12-\tf1{p_j})-\tf12,0\}
$$
is sufficient to ensure \eqref{ePCBB10}, as indicated by Theorem~\ref{thm:countexa}.

\end{rmk}

\begin{rmk}
Our result improves \cite[Theorem 2.1]{Kaur2020} since we allow $p_1$ or $p_2$ being less than 2 (this is only considered in the case of $n=1$ in \cite{Kaur2020}), and also $\al(p_1,p_2)< \al_{p_c}(p_1,p_2)$ when $2<p_1, p_2<p_c$. 
\end{rmk}

To obtain this result, we shall establish weighted and unweighted estimates of $B^{2,\al}_*$.
Recalling \eqref{eLBR7}, as the boundedness of $\tilde G^{1,\de_1}$ is clear, it remains to study the  boundedness of $B^{\be_1,\de_1}_{j,*}$.

\begin{prop}
	(i) Let $1<p\le 2$ and $n\ge 2$. Then, for $\be_1>n(\tf1p-\tf12)+\tf12$ and any $\ep>0$, there exists some $C_\ep>0$ such that
	\begin{equation}\label{ePCBB5}
		\|B^{\be_1,\de_1}_{j,*}(f)\|_{L^p(\rn)}\le C_\ep2^{-j(\al-\tf14-n(\tf1p-\tf12)-\ep)}\|f\|_{L^p(\rn)}.
	\end{equation}
	
	(ii) Let $0\le \ga<n$, where $n\ge 2$. For $\be_1>\max\{\tf{\ga-1}2,0\}+\tf12$ and any $\ep>0$, there exists some $C_\ep>0$ such that
	\begin{equation}\label{ePCBB6}
	\|B^{\be_1,\de_1}_{j,*}(f)\|_{L^2(|x|^{-\ga})}\le C_\ep2^{-j(\al-\tf14-\max\{\tf{\ga-1}2,0\}-\ep)}\|f\|_{L^2(|x|^{-\ga})}.
	\end{equation}
	
\end{prop}

\begin{rmk}
The operator $B^{\be_1,\de_1}_{j,*}$ was introduced by \cite{Kaur2020}, where they obtained its $L^p$-boundedness, when $p\ge 2$, by assuming the boundedness of $G^\psi_{j}$. It seems that \eqref{ePCBB5} is the first result  when $p<2$.
\end{rmk}

\bpf

(i)
We follow the argument in \cite[Section 5]{Kaur2020}. 
By definition, we control $B^{\be_1,\de_1}_{j,*}(f)(x)$ by
\begin{align} &\sup_{R>0}\Big(\int_0^{2^{-(j+1+\ep_0)/2}}|B_{j,t,R}^{\be_1}(f)(x)t^{2\de_1+1}|^2 dt\Big)^{1/2}\label{ePCBB2}\\
+&\sup_{R>0}\Big(\int_{2^{-(j+1+\ep_0)/2}}^{2^{-(j-1)/2}}|B_{j,t,R}^{\be_1}(f)(x)t^{2\de_1+1}|^2 dt\Big)^{1/2}.\label{ePCBB3}
\end{align}

To handle the first term \eqref{ePCBB2}, we have a pointwise control. Let us sketch it when $\be_1\in(0,1)$, while other cases can be dealt with similarly. Let $\rho=1-\be_1$ and $\psi^\kappa(x):=x^{-\kappa-\rho}\psi(x)$, then
$$
|B_{j,t,R}^{\be_1}(f)(x)|
\le 2^{j\rho}\sum_{\kappa\ge 0}\tf{\Ga(\rho+\kappa)}{\Ga(\rho)\kappa!}(2^jt^2)^\kappa|B^{\psi^\kappa}_{2^{-j},R}(f)(x)|,
$$
which implies further that
\begin{align}
	&\sup_{R>0}\Big(\int_0^{2^{-(j+1+\ep_0)/2}}|B_{j,t,R}^{\be_1}(f)(x)t^{2\de_1+1}|^2 dt\Big)^{1/2}\notag\\
	\le&2^{j\rho}\sum_{\kappa\ge 0}\tf{\Ga(\rho+\kappa)}{\Ga(\rho)\kappa!}2^{j\kappa}\sup_{R>0}\Big(\int_0^{2^{-(j+1+\ep_0)/2}}|B^{\psi^\kappa}_{2^{-j},R}(f)(x)t^{2\kappa}t^{2\de_1+1}|^2 dt\Big)^{1/2}\notag\\
	\le&2^{j\rho}\sum_{\kappa\ge 0}\tf{\Ga(\rho+\kappa)}{\Ga(\rho)\kappa!}2^{-\kappa(1+\varepsilon_0)}2^{-j(\de_1+\tf34)}\sup_{R>0}|B^{\psi^\kappa}_{2^{-j},R}(f)(x)|.\label{ePCBB1}
\end{align}
As 
$$
\sup_{x\in[\tf12,1],\ 0\le l\le N }\Big|\f{d^l\psi^\kappa}{dx^l}\Big|\le C(\rho)2^{N+\kappa}k^{N},
$$
it follows from 
\eqref{eLBR14} that, for $1<p\le 2$,
$$
\|\sup_{R>0}|B^{\psi^\kappa}_{2^{-j},R}(f)\|_{L^p(\rn)}
\le C 2^{N+\kappa}\kappa^{N} 2^{j((n-1)(\tf1p-\tf12)+\ep)}\|f\|_{L^p(\rn)}.
$$
This combined with \eqref{ePCBB1} yields 
\begin{align*}
\Big\|\sup_{R>0}\Big(\int_0^{2^{-(j+1+\ep_0)/2}}&|B_{j,t,R}^{\be_1}(f)(x)t^{2\de_1+1}|^2 dt\Big)^{1/2}\Big\|_{L^p(\rn)}\\
\le &C2^{j\rho}\sum_{\kappa\ge 0}\tf{\Ga(\rho+\kappa)}{\Ga(\rho)\kappa!}2^{-j(\de_1+\tf34)}\kappa^{N}2^{-\kappa\varepsilon_0} 2^{j((n-1)(\tf1p-\tf12)+\ep)}\|f\|_{L^p(\rn)}\\
\le &C 2^{-j(\al-\tf14-(n-1)(\tf1p-\tf12)-\ep)}\|f\|_{L^p(\rn)}
\end{align*}
as the summation over $\kappa$ is convergent. 

We discuss the second term \eqref{ePCBB3} below.
We sketch the proof in \cite{Kaur2020}. By change of variables, we control \eqref{ePCBB3} by 
\begin{align*}
	&2^{-j(\de_1+\tf14)}\sup_{R>0}\Big(\int_{s_1}^{s_2}|[(s^2-\tf{|\xi|^2}{R^2})^{\be_1-1}_+\psi(2^j(1-\tf{|\xi|^2}{R^2}))\wh f(\xi)]^\vee(x)|^2ds\Big)^{1/2}\\
	\le &2^{-j(\de_1+\tf14)}\sum_{\kappa\ge j-2} 2^{-\kappa(\be_1-1)}\sup_{R>0}\Big(\int_{s_1}^{s_2}|[\psi(2^\kappa(1-\tf{|\xi|^2}{R^2s^2}))\psi(2^j(1-\tf{|\xi|^2}{R^2}))\wh f(\xi)]^\vee(x)|^2ds\Big)^{1/2},
\end{align*}
where $s_1^2=1-2^{-(j-1)}$, and $s_2^2=1-2^{-(j+1+\ep_0)}$.
We can decompose the interval $[2^{-\kappa-1},2^{-\kappa+1}]$ into subintervals $I_m=c_m+O(2^{-\kappa(1+\ep)})$ of length $ 2^{-\kappa(1+\ep)}$, where $c_m=1-2^{-\kappa+1}+2^{-\kappa(1+\ep)}m\sim 1$, $0\le m\le C2^{\kappa\ep}$. Let $\vp$ be a compactedly supported smooth function such that
$$
\chi_{[2^{-\kappa-1},2^{-\kappa+1}]}=\chi_{[2^{-\kappa-1},2^{-\kappa+1}]}\sum_{m=0}^{C2^{\kappa\ep}}\vp(2^{\kappa(1+\ep)}(c_m-|\xi|^2)).
$$
It follows from Taylor's expansion that the main term of 
$$
\psi(2^\kappa(1-|\xi|^2))\psi(2^j(1-s^2|\xi|^2))\vp(2^{\kappa(1+\ep)}(c_m-|\xi|^2))
$$
is 
\begin{equation}\label{ePCBB4}
	\psi(2^\kappa(1-c_m))\psi(2^j(1-s^2 c_m))\vp(2^{\kappa(1+\ep)}(c_m-|\xi|^2)).
\end{equation}
Applying \eqref{ePCBB4}, by a change of variables, we essentially have the control 
\begin{align*}
\eqref{ePCBB3}\les &2^{-j(\de_1+\tf14)}\sum_{\kappa\ge j-2} 2^{-\kappa(\be_1-1)}2^{\kappa\ep}\sup_{R>0}\Big(\int_{s_1}^{s_2} \Big|B^\vp_{(2^{\kappa(1+\ep)}c_m)^{-1},sR}(f)(x)\Big|^2 ds\Big)^{1/2}\\
\les&2^{-j(\de_1+\tf14)}\sum_{\kappa\ge j-2} 2^{-\kappa(\be_1-1)}2^{\kappa\ep}
\tilde G^\vp_{\kappa(1+\ep)}(f)(x), 
\end{align*}
where we use $s\sim 1$.

Applying \eqref{eLBR19}, we have, for $1<p\le 2$,
$$
\|\tilde G^\vp_{\kappa(1+\ep)}(f)\|_{L^p(\rn)}
\le C2^{\kappa(1+\ep)[n(\tf1p-\tf12)-\tf12]}\|f\|_{L^p(\rn)}.
$$
Therefore the contribution of \eqref{ePCBB3} is bounded by 
$$
2^{-j(\al-\tf14-n(\tf1p-\tf12)-\ep)}\|f\|_{L^p(\rn)}
$$
when $1<p\le 2$ and $\be_1>n(\tf1p-\tf12)+\tf12$.
This and the estimate of \eqref{ePCBB2} complete the proof of \eqref{ePCBB5}.

(ii)
As in the estimates of the $L^p$-norm, we need to estimate the contribution of \eqref{ePCBB2} and  \eqref{ePCBB3} separately. 
It follows from the pointwise control \eqref{ePCBB1} and Theorem~\ref{tLBR5} that
\begin{align*}
	\Big\|\sup_{R>0}\Big(\int_0^{2^{-(j+1+\ep_0)/2}}&|B_{j,t,R}^{\be_1}(f)(x)t^{2\de_1+1}|^2 dt\Big)^{1/2}\Big\|_{L^2(|x|^{-\ga})}\\
	\les&
	2^{j\rho}\sum_{\kappa\ge 0}\tf{\Ga(\rho+\kappa)}{\Ga(\rho)\kappa!}2^{-j(\de_1+\tf34)}\kappa^{N}2^{-\kappa\varepsilon_0} 2^{j/2}A_j(\ga)\|f\|_{L^2(|x|^{-\ga})}\\
	\les& 
	2^{-j(\al-\tf34-\max\{\tf\ga2-1,-\tf12\}-\ep)}\|f\|_{L^2(|x|^{-\ga})}.
\end{align*}
To handle \eqref{ePCBB3},
we apply \eqref{ePCBB4} to obtain an essential control 
$$
\eqref{ePCBB3}\le 2^{-j(\de_1+\tf14)}\sum_{k\ge j-2} 2^{-\kappa(\be_1-1)}2^{\kappa\ep}\Big(\int_0^\nf \Big|B^\vp_{(2^{\kappa(1+\ep)}c_m)^{-1},s}(f)(x)\Big|^2\f{ds}s\Big)^{1/2}.
$$
From Theorem~\ref{tLBR5} we know 
$$
\left\|\Big(\int_0^\nf \Big|B^\vp_{(2^{\kappa(1+\ep)}c_m)^{-1},s}(f)(x)\Big|^2\f{ds}s\Big)^{1/2}\right\|_{L^2(|x|^{-\ga})}
\le CA_{\kappa(1+\ep)}(\ga)\|f\|_{L^2(|x|^{-\ga})}.
$$
As a result, the contribution of of \eqref{ePCBB3}, when $\be_1>\max\{\tf{\ga-1}2,0\}+\tf12$, is bounded by 
$$
2^{-j(\al+\tf14-\max\{\tf{\ga-1}2,0\}-\tf12-\ep)}\|f\|_{L^2(|x|^{-\ga})}.
$$

We complete the proof of (ii).

\epf

We are ready to prove Theorem~\ref{tPCBB1}, which is contained in the following three results.

\begin{prop}\label{tPCBB2}

	(i) Let $n\ge 2$ and $2\le p_1,p_2<\nf$. Suppose $\al>\max\{n(\tf12-\tf1{p_1})-\tf12,0\}+\max\{n(\tf12-\tf1{p_2})-\tf12,0\}$. Then, for any $f_1\in L^{p_1}(\rn)$ and $f_2\in L^{p_2}(\rn)$, we have
	$$
	\lim_{R\to\nf} B^{2,\al}_R(f_1,f_2)(x)= f_1(x)f_2(x)\quad a.e.
	$$
	
	(ii) Let $0\le \ga_1,\ga_2<n$, where $n\ge 2$, and $\ga=\tf{\ga_1+\ga_2}2$.
	Suppose 
	\begin{equation}\label{ePCBB7}
		\al>\max\{\tf{\ga_1-1} 2,0\}+\max\{\tf{\ga_2-1} 2,0\}.
	\end{equation} 
Then
	$$
	\|B^{2,\al}_*(f_1,f_2)\|_{L^1(|x|^{-\ga})}\le C\|f_1\|_{L^2(|x|^{-\ga_1})}\|f_2\|_{L^2(|x|^{-\ga_2})}.
	$$

\end{prop}

\bpf

(i)  By Stein's maximal principle in the bilinear setting (see, for instance, \cite[Proposition 4.4]{Grafakos2018a}), (i) follows from (ii) and that 
\begin{equation}\label{ePCBB8}
L^q(\rn)\sset L^2(\rn)+L^2(|x|^{-\om})
\end{equation}
for $\om>n(1-\tf 2q)$ when $q\ge 2$.

We prove (ii) below.

(ii) For $j\ge 1$, it follows from \eqref{eLBR7}, \eqref{eLBR8}, and \eqref{ePCBB6} that, for $\al$ satisfying \eqref{ePCBB7},
\begin{align*}
	\|M_j(f_1,f_2)\|_{L^1(|x|^{-\ga})}
	&\les  2^{-j/4}\|\tilde G^{\de_1}(f_1)\|_{L^2(|x|^{-\ga_1})}\|B^{\be_1,\de_1}_{j,*}(f_2)\|_{L^2(|x|^{-\ga_2})}\\
	&\les 2^{-j(\al-\max\{\tf{\ga_2-1}2, 0\}-\ep)}\|f_1\|_{L^2(|x|^{-\ga_1})}\|f_2\|_{L^2(|x|^{-\ga_2})},
\end{align*}
where $\be_1$ and $\de_1$ are chosen so that 
$\be_1+\de_1=\al$, $\be_1>\max\{\tf{\ga_2-1} 2,0\}+\frac 12$, and $\de_1>\max\{\tf{\ga_1-1} 2,0\}-\frac 12$.
As $\al -\max\{ \tf{\ga_2-1}2, 0\}-\ep>0$, we obtain
$$
\|\sum_{j\ge1}M_j(f_1,f_2)\|_{L^1(|x|^{-\ga})}\les \|f_1\|_{L^2(|x|^{-\ga_1})}\|f_2\|_{L^2(|x|^{-\ga_2})}.
$$

It remains to handle $M_0$ associated with $m^\al_0$ in \eqref{eLBR5}. We can actually decompose
$$
m^\al_0=\sum_{j\ge0} m_0^j,
$$
where $m_0^j(\xi)=m^\al_0(\xi)\psi(2^j(1-|\xi_1|^2))$, and
$m_0^0=m^\al_0-\sum_{j\ge 1}m_0^j$. The study of $m_0^j$ is similar to that of $m^\al_j$ with the roles of $f_1$ and $f_2$ switched, and $m_0^0\in C_c^\nf(\mathbb R^{2n})$,
therefore we have
$$
\|M_0(f_1,f_2)\|_{L^1(|x|^{-\ga})}\les \|f_1\|_{L^2(|x|^{-\ga_1})}\|f_2\|_{L^2(|x|^{-\ga_2})}
$$
when $	\al>\max\{\tf{\ga_1-1} 2,0\}+\max\{\tf{\ga_2-1} 2,0\}$.

We complete the proof of (ii).

\epf

\begin{prop}\label{prop:12le2}
	Let $n\ge 2$ and $1<p_1, p_2\le 2$, and $\al>n(\tf1p-1)$.  Then, for any $f_1\in L^{p_1}(\rn)$ and $f_2\in L^{p_2}(\rn)$, we have
	$$
	\lim_{R\to\nf} B^{2,\al}_R(f_1,f_2)(x)= f_1(x)f_2(x)\quad a.e.
	$$
\end{prop}

\bpf

The proof is similar to Proposition~\ref{tPCBB2}.

For $j\ge 1$, we apply \eqref{eLBR7},  \eqref{eLBR20}, and \eqref{ePCBB5} to obtain that
\begin{align*}
	\|M_j(f_1,f_2)\|_{L^p(\rn)}
	\les & 2^{-j/4}\|\tilde G^{\de_1}(f_1)\|_{L^{p_1}(\rn)}\|B^{\be_1,\de_1}_{j,*}(f_2)\|_{L^{p_2}(\rn)}\\
	\les& 2^{-j(\al-n(\tf1{p_2}-\tf12)-\ep)}\|f_1\|_{L^{p_1}(\rn)}\|f_2\|_{L^{p_2}(\rn)}
\end{align*}
when $\al=\de_1+\be_1>n(\tf1{p_1}-\tf12)-\tf12+n(\tf1{p_2}-\tf12)+\tf12=n(\tf1p-1)>n(\tf1{p_2}-\tf12)+\ep$.
Summing over $j$ we obtain 
$$
	\|B^{2,\al}_*(f_1,f_2)\|_{L^{p}(\rn)}\le C\|f_1\|_{L^{p_1}(\rn)}\|f_2\|_{L^{p_2}(\rn)},
$$
which implies the conclusion.

\epf

\begin{prop}
	Let $1<p_1\le 2\le p_2<\nf$, and $\al>n(\tf1{p_1}-\tf12)+\max\{0,n(\tf12-\tf1{p_2})-\tf12\}$.
	Then, for any $f_1\in L^{p_1}(\rn)$ and $f_2\in L^{p_2}(\rn)$, we have
	\begin{equation}\label{ePCBB9}
	\lim_{R\to\nf} B^{2,\al}_R(f_1,f_2)(x)\to f_1(x)f_2(x)\quad a.e.
	\end{equation}
\end{prop}

\bpf
Let $\tf1q=\tf1{p_1}+\tf12$.   Take $\gamma_2> n(1-\frac 2{p_2})$ so that 
$
\al>n(\tf1{p_1}-\tf12)+\max\{\tf{\ga_2-1}2,0\}
$ and we get the splitting $f_2=f_2^1+f_2^2$, where $f_2^1 \in L^2$ and $f_2^2 \in L^2(|x|^{-\gamma_2})$. Since 
\[
\al>n(\tf1{p_1}-\tf12)+\max\{0,n(\tf12-\tf1{p_2})-\tf12\}\ge n(\tf1{p_1}-\tf12)+n(\frac 12-\frac 12),
\]by Proposition \ref{prop:12le2},
we already know that 
\[
\lim_{R\to\nf} B^{2,\al}_R(f_1,f_2^1)(x)\to f_1(x)f_2^1(x)\quad a.e.
\]
It remains to prove 
\[
\lim_{R\to\nf} B^{2,\al}_R(f_1,f_2^2)(x)\to f_1(x)f_2^2(x)\quad a.e.
\]

To handle the case $j\ge 1$, we apply \eqref{eLBR7}, \eqref{eLBR20}, and \eqref{ePCBB6} to obtain
$$
\|M_j(f_1,f_2^2)\|_{L^q(|x|^{-\ga_2q/2})}
	\les 2^{-j(\al-\max\{\tf{\ga_2-1}2,0\}-\ep)}\|f_1\|_{L^{p_1}(\rn)}\|f_2^2\|_{L^{2}(|x|^{-\ga_2})}
$$
since $\al>n(\tf1{p_1}-\tf12)+\max\{\tf{\ga_2-1}2,0\}$. 
It follows that 
$$
\|\sum_{j\ge 1}M_j(f_1,f_2^2)\|_{L^q(|x|^{-\ga_2q/2})}
	\les \|f_1\|_{L^{p_1}(\rn)}\|f_2^2\|_{L^{2}(|x|^{-\ga_2})}.
$$
By a similar argument, one can obtain from \eqref{ePCBB5} and \eqref{eLBR8} that
$$
\|M_0(f_1,f_2^2)\|_{L^q(|x|^{-\ga_2q/2})}
	\les \|f_1\|_{L^{p_1}(\rn)}\|f_2^2\|_{L^{2}(|x|^{-\ga_2})}.
$$
We conclude from \eqref{ePCBB8} and Stein's maximal principle that \eqref{ePCBB9} holds under our assumptions.

\epf


\section{Multilinear Bochner-Riesz means}

In the study of bilinear Bochner-Riesz operators, the square functions play an important role. 
In this section, we establish the relation between $k$-linear Bochner-Riesz means and $(k-1)$-linear square functions. As applications of this relation, we obtain the boundedness and pointwise convergence of multilinear Bochner-Riesz means.

The study of $B_*^{k,\al}$ is closely related to the associated  multilinear square function
\begin{equation*} 
\tilde G^{k,\al}(f_1,\dots, f_k)(x)=\sup_{R>0}\Big(\f 1R\int_0^R |B^{k,\al}_t(f_1,\dots, f_k)(x)|^2 dt\Big)^{1/2}.
\end{equation*}
To explain this, we decompose 
\begin{equation*} 
m^\al=\sum_{j\ge 0}m^\al_j,
\end{equation*}
where $m^\al_j(\xi)=(1-|\xi|^2)^{\al}_+ \psi(2^j(1-|\xi_k|^2))$ for $j\ge 1$ and $m^\al_0(\xi)=m^\al(\xi)-\sum_{j\ge 1}m^\al_j(\xi)$. The maximal function related to $m^\al_j$ is defined by 
$$
M_j(f_1,\dots,f_k)(x)=\sup_{R>0}\Big|\int_{\bbr^{nk}}m^\al_j(R^{-1}\xi)\prod_{\iota=1}^k \wh f_\iota(\xi_\iota)e^{2\pi ix\cdot(\xi_1+\cdots+\xi_k)} d\xi\Big|.
$$

Similar to \eqref{eLBR21}, we have the following result.

\begin{prop}
Let $j\ge 1$. Then
\begin{align}
&\Big|\int_{\bbr^{nk}}m^\al_j(R^{-1}\xi)\prod_{\iota=1}^k \wh f_\iota(\xi_\iota)e^{2\pi ix\cdot(\xi_1+\cdots+\xi_k)} d\xi\Big|\notag\\
\les& 2^{-j/4} \tilde G^{k-1,\de}(f_1,\dots, f_{k-1})(x) \Big(\int_0^{2^{-(j-1)/2}}|B^{\be}_{j,t,R}(f_k)(x)t^{2\de+1}|^2dt\Big)^{1/2}.\label{eMBR1}
\end{align}
As a result,
\begin{equation}\label{eMBR8}
M_j(f_1,\dots, f_k)(x)
\les 2^{-j/4}\tilde G^{k-1,\de}(f_1,\dots, f_{k-1})(x)
B_{j,*}^{\be, \delta}(f_k)(x).
\end{equation}

\end{prop}

Since  $B_{j,*}^{\be, \delta}$  defined in \eqref{eLBR6} is well understood, we will focus on the study of $\tilde G^{{k-1},\de}$. 
For $j\ge 0$, we define
$$
\tilde G^{k,\al}_j(f_1,\dots, f_k)(x)=\sup_{R>0}\Big(\f 1R\int_0^R \Big|\int_{\bbr^{nk}}m^\al_j(t^{-1}\xi)\prod_{\iota=1}^k \wh f_\iota(\xi_\iota)e^{2\pi ix\cdot(\xi_1+\cdots+\xi_k)} d\xi\Big|^2 dt\Big)^{1/2},
$$
then 
$\tilde G^{k,\al}\les\sum_{j\ge 0}\tilde G^{k,\al}_j$.

For $j\ge 1$, we apply \eqref{eMBR1} to obtain 
\begin{align}
\tilde G^{k,\al}_j(f_1,\dots, f_k)(x)
&\les  2^{-j/4}\tilde G^{k-1,\de} (f_1,\dots, f_{k-1})(x)\notag\\
&\q \times\sup_{L>0}\Big(\tf 1L\int_0^L \int_0^{2^{-(j-1)/2}}|B^{\be}_{j,t,R}(f_k)(x)t^{2\de+1}|^2dt dR\Big)^{1/2}\notag\\
&  \le 2^{-j/4}\tilde G^{k-1,\de} (f_1,\dots, f_{k-1})(x)B_{j,*}^{\be,\de}(f_k)(x).\label{eMBR5}
\end{align}
Iterating the above estimate we get 
\begin{align*}
\tilde G^{k,\al}_j(f_1,\dots, f_k)(x)
&\q\les 2^{-j/4}\tilde G^{k-1,\de_k}(f_1,\dots, f_{k-1})(x)
B_{j,*}^{\be_k, \delta_k}(f_k)(x)\\
&\q\lesssim  \sum_{j_{k-1}\ge 0}\cdots \sum_{j_2\ge 0} 2^{-(j_2+\cdots+ j_{k-1}+j_k)/4}\tilde G^{1,\be_1}(f_1) \prod_{i=2}^{k}B_{j_\iota,*}^{\be_\iota, \delta_\iota}(f_\iota)(x),
\end{align*}
where $j_k:=j$ and
\begin{equation*} 
\alpha=\delta_k+\beta_k= \sum_{\iota=1}^{k} \beta_\iota,\qquad \delta_{\iota+1}=\delta_\iota+\beta_\iota,\q \de_2=\be_1.
\end{equation*}
Note that $M_j$ satisfies the same estimate.

We denote $I_k=\{1,\dots,k\}$, and  assume that $I\sset I_k$.

\begin{thm}\label{tMBR2}
Let $I\sset I_k$, $\ga_\iota\in(0,n)$ for $\iota\in I$, and
	 $p_\iota\in (1,2]$ for $\iota\in I_k\setminus I$.
	Suppose 
	 $\tf1p=\sum_{\iota\in I}\tf12+\sum_{\iota\in I_k\setminus I}\tf1{p_\iota}$,
	 $\ga=\tf p{2}\sum_{\iota\in I}\ga_\iota$, and
	$$
	\al>\sum_{\iota\in I}\max\{\tf{\ga_\iota-1}2,0\}+\sum_{\iota\in I_k\setminus I}n(\tf1{p_\iota}-\tf12)+\tf{k-2}2.
	$$
	Then 
	\begin{equation}\label{eMBR7}
	\|\tilde G^{k,\al}(f_1,\dots, f_k)\|_{L^p(|x|^{-\ga})}
	\les \prod_{\iota\in I}\|f_\iota\|_{L^2(|x|^{-\ga_\iota})}
	\prod_{\iota\in I_k\setminus I} \|f_\iota\|_{L^{p_\iota}(\rn)}.
	\end{equation}
\end{thm}

\bpf
We prove \eqref{eMBR7} inductively.

For $k=1$, \eqref{eMBR7} follows from \eqref{eLBR20} and \eqref{eLBR8}.

Suppose \eqref{eMBR7} holds for $k-1$. Moreover, we may assume that $k\notin I$ while the case $k\in I$ can be handled similarly. 
Applying \eqref{eMBR5}, we see that, for $j\ge 1$, 
$$
\|\tilde G^{k,\al}_j(f_1,\dots, f_k)\|_{L^p(|x|^{-\ga})}
\les 2^{-j/4} \|\tilde G^{k-1,\de}(f_1,\dots, f_{k-1})\|_{L^{q}(|x|^{-\ga_q})}\|B^{\be,\de}_{j,*}(f_k)\|_{L^{p_k}(\rn)},
$$
where 
$\tf1q=\sum_{\iota\in I}\tf12+\sum_{\iota\in I_{k-1}\setminus I}\tf1{p_\iota}=\tf1p-\tf1{p_k}$, $\ga_q=\tf{\ga q}p=\tf q{2}\sum_{\iota\in I}\ga_\iota$,
$\be>n(\tf1p-\tf12)+\tf12$, and $\de=\al-\be>\sum_{\iota\in I}\max\{\tf{\ga_\iota-1}2,0\}+\sum_{\iota\in I_{k-1}\setminus I}n(\tf1{p_\iota}-\tf12)+\tf{k-3}2$. By the inductive hypothesis, we have
$$
 \|\tilde G^{k-1,\de}(f_1,\dots, f_{k-1})\|_{L^{q}(|x|^{-\ga_q})}
 \les \prod_{\iota\in I}\|f_\iota\|_{L^2(|x|^{-\ga_\iota})}
	\prod_{\iota\in I_{k-1}\setminus I} \|f_\iota\|_{L^{p_\iota}(\rn)}.
$$
This combined with \eqref{ePCBB5} yields the bound
$$
\|\tilde G^{k,\al}_j(f_1,\dots, f_k)\|_{L^p(|x|^{-\ga})}
\les 2^{-j(\al-n(\tf1{p_k}-\tf12)-\ep)} \prod_{\iota\in I}\|f_\iota\|_{L^2(|x|^{-\ga_\iota})}
	\prod_{\iota\in I_k\setminus I} \|f_\iota\|_{L^{p_\iota}(\rn)}.
$$
Summing over $j\ge 1$, we obtain 
$$
\|\sum_{j\ge 1}\tilde G^{k,\al}_j(f_1,\dots, f_k)\|_{L^p(|x|^{-\ga})}
\les  \prod_{\iota\in I}\|f_\iota\|_{L^2(|x|^{-\ga_\iota})}
	\prod_{\iota\in I_k\setminus I} \|f_\iota\|_{L^{p_\iota}(\rn)}.
$$

To handle $\tilde G^{k,\al}_0$, we need to decompose the multiplier $m^\al_0$ further in $\xi_{k-1}$, $\xi_{k-2},\dots,\ \xi_1$ successively, as we did in Proposition~\ref{tPCBB2}. We leave the details to the reader.

\epf

\bpf[Proof of Theorem~\ref{tMBR6}]

By comparing \eqref{eMBR8} and \eqref{eMBR5}, we see that the conclusion follows from
Theorem~\ref{tMBR2} by modifying the inductive argument in Theorem~\ref{tMBR2} .

\epf

The following result is essentially proved by \cite{Shuin2024}. We can prove it using the relation \eqref{eMBR5} and  Proposition~\ref{tLBR2}, which is similar to the proof of Theorem~\ref{tMBR2}. Compared with the proof in \cite{Shuin2024}, our approach avoids using the boundedness of the operators $A^{R,t}_{j,\be}$, which is the main obstacle to overcome in \cite{Choudhary2023}.

\begin{prop}[\cite{Shuin2024}]

Let $p_c>\tf{2n}{n-1}$, $p_\iota\in [2,\nf)$ for $\iota=1,\dots,k$.
Let $\al_{p_c}(p_1,\dots, p_k)=\sum_{\iota\in I}\al(p_\iota)+\sum_{\iota\in I_k\setminus I}\al_{p_c}(p_\iota)$, where  $I:=\{\iota\in I_k:\ p_\iota\ge p_c\}$
Suppose that for any $\ep>0$ there exists some $C_\ep>0$ such that
\begin{equation*}
\Big\|\Big(\int_1^2\Big|[\psi(2^j(t^2-|\xi|^2))\wh f(\xi)]^\vee(x)\Big|^2dt\Big)^{1/2}\Big\|_{L^p(\rn)}\le C_\ep 2^{-j(\tf12-\al(p)-\ep)}\|f\|_{L^p(\rn)}
\end{equation*}
holds for all $p>p_c$.
 Then
$$
\|B_*^{k,\al}(f_1,f_2,\dots, f_k)\|_{L^p}\le C\prod_{\iota\in I}\|f_\iota\|_{L^{p_\iota}}
$$
whenever $\al>\al_{p_c}(p_1,\dots, p_k)$.

\end{prop}

We remark that, as in the bilinear case, 
$$
\al_{p_c}(p_1,\dots, p_k)>\al(p_1,\dots, p_k)
$$
when $p_\iota\in (2,p_c)$ for some $\iota$.

\section{A counterexample}

We prove Theorem~\ref{thm:countexa} in this section.

Noting that
$$\int_{\bbr^m}(1-|\xi|^2)_+^\alpha e^{2\pi ix\cdot \xi}\;d\xi=C\frac{J_{\frac{m}2+\alpha}(2\pi|x|)}{|x|^{\frac{m}2+\alpha}},$$
 we get
\begin{align*}\nonumber
&B^{k,\al}_R(f_1,\dots,f_k)(x)\\
=&\int_{\bbr^{nk}}\left\{\int_{\bbr^{nk}}\left(1-\frac{|\xi|^2}{R^2}\right)_+^\alpha e^{2\pi i(x-z_1,\cdots,x-z_k)\cdot \xi}\;d\xi\right\}\prod_{\iota=1}^kf_\iota(z_\iota)\;dz_1\cdots dz_k\\ 
=&\int_{\bbr^{nk}} K_{k,R}^\alpha(x-z_1,\cdots,x-z_k)\prod_{\iota=1}^kf_\iota(z_\iota)\;dz_1\cdots dz_k, \nonumber
\end{align*}
where 
$$
K_{k,R}^\alpha(z_1,\cdots,z_k)=R^{nk}\frac{J_{\frac{nk}2+\alpha}(2\pi R|(z_1,\cdots,z_k)|)}{(R|(z_1,\cdots,z_k)|)^{\frac{nk}2+\alpha}}.
$$

Let $M$ be a large number and $\epsilon$ be a small number. We define $\varphi(x)$ to be a smooth function by
$$\varphi(x)=\varphi_{\epsilon,M}(x)=\psi(\epsilon^{-1}|x'|)\psi(\epsilon^{-1}M^{-\frac12}x_n),\; x=(x',x_n)\in\bbr^{n-1}\times\bbr,$$
where $\psi$ is a smooth bump function supported in the interval $[-1,1]$. We also denote the function 
$$f_{\epsilon,M}(x)=e^{2\pi ix_n}\varphi(x),$$
and the set $S_M$ 
\begin{equation*} 
S_M:=\{x\in\bbr^n:\;M\leq|x'|\leq 2M,\;M\leq x_n\leq 2M\}.
\end{equation*}
Then, it is easy to check that
\begin{equation}\label{equ:fmlp}
\|f_{\epsilon,M}\|_{L^p(\bbr^n)}\simeq (\epsilon^n M^\frac12)^\frac1p,\quad \text{and}\quad |S_M|\simeq M^n.
\end{equation}

For $x\in S_M$, we claim that for  $R=\frac{\sqrt k|x|}{x_n}\simeq 1$, 
\begin{equation}\label{equ:Bkalpclaim}
\big| B^{k,\al}_R(f_{\epsilon,M},\dots,f_{\epsilon,M})(x)\big|\geq CM^{-\frac{nk+1}2-\alpha}(\epsilon^n M^\frac12)^k. 
\end{equation}
With this claim in hand, we get
\begin{align}\label{equ:lowbound}
\big\|B^{k,\al}_*(f_{\epsilon,M},\dots, f_{\epsilon,M})(x)\big\|_{L^p(\bbr^n)}\geq&CM^{-\frac{nk+1}2-\alpha+\frac{k}2} |S_M|^\frac1p\\\nonumber
\geq& CM^{-\frac{nk+1}2-\alpha+\frac{k}2+\frac{n}p}.
\end{align}
On the other hand, by \eqref{equ:fmlp}, we have
$$\prod_{\iota=1}^k\|f_{\epsilon,M}\|_{L^{p_\iota}(\bbr^n)}\simeq  (\epsilon^n M^\frac12)^{\frac1{p_1}+\cdots\frac1{p_k}}.$$
This together with \eqref{equ:lowbound}   yields  that
\begin{align*}
M^{-\frac{nk+1}2-\alpha+\frac{k}2+\frac{n}p}\lesssim  M^{\frac12(\frac1{p_1}+\cdots\frac1{p_k})}.
\end{align*}
So we get the desired result \eqref{equ:condialp} by letting $M\to\infty$. It remains to verify \eqref{equ:Bkalpclaim}.

\bpf[Proof of the claim \eqref{equ:Bkalpclaim}]
By the definition of $S_M$ and the support property of $\varphi$, we obtain that
for $x\in S_M$ and $z_\iota\in\; {\rm supp} f_{\epsilon,M}$, 
\begin{equation}\label{equ:xzj}
|x-z_\iota|\simeq M\gg1.
\end{equation}
Combining this with the asymptotic behavior of Bessel function
$$J_m(r)=C(m)\frac1{\sqrt{r}}e^{ir}+\overline{C(m)}\frac1{\sqrt{r}}e^{-ir}+O(r^{-\frac32}),\quad \text{as}\quad r\to+\infty,$$
we deduce that for $x\in S_M,\;R=\frac{\sqrt k|x|}{x_n}\simeq 1$,
\begin{align*}
&B^{k,\al}_R(f_{\epsilon,M},\dots,f_{\epsilon,M})(x)\\
=&C_1R^{nk}\int_{\bbr^{nk}}\frac{e^{2\pi iR|(x-z_1,\cdots,x-z_k)|}e^{2\pi i(z_{1,n}+\cdots+z_{k,n})}}{(R|(x-z_1,\cdots,x-z_k)|)^{\frac{nk+1}2+\alpha}}
\prod_{\iota=1}^k\varphi(z_\iota)\;dz_1\cdots dz_k\\
&+C_2R^{nk}\int_{\bbr^{nk}}\frac{e^{-2\pi iR|(x-z_1,\cdots,x-z_k)|}e^{2\pi i(z_{1,n}+\cdots+z_{k,n})}}{(R|(x-z_1,\cdots,x-z_k)|)^{\frac{nk+1}2+\alpha}}
\prod_{\iota=1}^k\varphi(z_\iota)\;dz_1\cdots dz_k\\
&+C_3R^{nk}\int_{\bbr^{nk}}O\left(\frac{1}{(R|(x-z_1,\cdots,x-z_k)|)^{\frac{nk+3}2+\alpha}} \right)\prod_{\iota=1}^k\varphi(z_\iota)\;dz_1\cdots dz_k\\
=:&\ I_1+I_2+I_3.
\end{align*}

{\bf The estimate of the error term $I_3$:}  By the definition of $\varphi(x)$, \eqref{equ:xzj} and $R\simeq 1$, we get
\begin{align*}
|I_3|\lesssim& M^{-\frac{nk+3}2-\alpha}(\epsilon^n M^\frac12)^k\lesssim \epsilon^{nk} M^{-\frac{nk+3}2-\alpha+\frac{k}2}.
\end{align*}

{\bf The estimate of the error term $I_2$:}  By integration by part in the variable $z_{1,n}$, we derive
\begin{align*}
|I_2|\lesssim& M^{-\frac{nk+1}2-\alpha}(\epsilon^n M^\frac12)^k\big(M^{-1}+\epsilon^{-1}M^{-\frac12}\big).
\end{align*}

{\bf The estimate of the main term $I_1$:}  Noting that
$$\big| R|(x-z_1,\cdots, x-z_k)|+z_{1,n}+\cdots+z_{k,n}-R|(x,\cdots,x)|\big|\leq C\epsilon$$
we get
$$|I_1|\geq CM^{-\frac{nk+1}2-\alpha}(\epsilon^n M^\frac12)^k.$$

Hence, we obtain the claim \eqref{equ:Bkalpclaim}.

\epf

\end{document}